\documentclass[12pt,a4paper,fleqn]{article}
\usepackage{a4wide,amsfonts,amsmath,latexsym,amssymb,euscript,graphicx,units,mathrsfs}

\usepackage{graphicx}
\usepackage{color}
\usepackage{amssymb}
\usepackage{amssymb}
\usepackage[T1]{fontenc}
\usepackage{latexsym}
\usepackage{xypic}
\usepackage{eufrak}
\usepackage{euscript}
\usepackage{amsfonts,amsmath}
\usepackage{verbatim}
\usepackage{fancyhdr}
\usepackage[english]{babel}
\usepackage{mathrsfs}
\usepackage{units}

\newtheorem{prop}{Proposition}[section]

\newtheorem{lemme}[prop]{Lemma}
\newtheorem{rem}[prop]{Remark}
\newtheorem{thm}[prop]{Theorem}
\newtheorem{defi}[prop]{Definition}

\renewcommand{\geq}{\geqslant}
\def\leq{\leqslant}

\newcommand{\Z}{\mathbb{Z}}
\newcommand{\R}{\mathbb{R}}

\def\HH{\EuFrak H}

\def\e{\varepsilon}

\def\1{{\mathbf{1}}}

\def\1{{\mathbf{1}}}
\def\0.5{{\frac{1}{2}}}

\newcommand{\fin}
{ \vspace{-0.6cm}
\begin{flushright}
\mbox{$\Box$}
\end{flushright}
\noindent }

\newcommand{\qed}{\nopagebreak\hspace*{\fill}
{\vrule width6pt height6ptdepth0pt}\par}

\begin{document}

\begin{center}
{\large{\bf 
A change of variable formula for the 2D fractional Brownian motion
of Hurst index bigger or equal to $1/4$
}}\\~\\
by Ivan Nourdin\footnote{Laboratoire 
de Probabilit{\'e}s et Mod{\`e}les Al{\'e}atoires, Universit{\'e} Pierre et Marie Curie,
Bo{\^\i}te courrier 188, 4 Place Jussieu, 75252 Paris Cedex 5, France,
{\tt ivan.nourdin@upmc.fr}}\\
{\it Universit\'e Paris VI}\\~\\
\small{{\it This version}: October 2, 2008}\\
\end{center}

{\small
\noindent
{\bf Abstract:} We prove a change of variable formula for the 2D fractional Brownian motion 
of index $H$ bigger or equal to $1/4$. For $H$ strictly bigger than $1/4$, our formula coincides with that obtained by using the rough paths theory. 
For $H=1/4$ (the more interesting case), there is an additional term that is a classical Wiener integral 
against an independent standard Brownian motion.\\



\section{Introduction and main result}
In \cite{coutin-qian}, Coutin and Qian have shown that
the rough paths theory of Lyons \cite{lyons} 
can be applied to the 2D fractional Brownian motion
$B=(B^{(1)},B^{(2)})$ under the condition that its Hurst parameter $H$ (supposed 
to be the same for the two components)
is {\sl strictly} bigger than $1/4$.
Since this seminal work, several authors 
have recovered this fact by using different routes
(see e.g. Feyel and de La Pradelle \cite{feyel-pradelle},
Friz and Victoir \cite{friz-victoir} or Unterberger \cite{unterberger} to cite but a few).
On the other hand, it is still an open  problem
to bypass this restriction  on $H$.

Rough paths theory is purely deterministic in essence. Actually, its random aspect 
comes only when it is applied to a single path of a given {\sl stochastic} process (like a 
Brownian motion, a fractional 
Brownian motion, etc.). In particular, {\sl it does not allow
to produce a new alea}. As such, the second point of Theorem \ref{main-thm} just below 
shows, in a sense, that it seems difficult to reach the case $H=1/4$ by using exclusively 
the tools of rough paths theory.

Before stating our main result, we need some preliminaries.
Let $W$ be a standard (1D) Brownian motion, independent of $B$.
We assume that $B$ and $W$ are defined on the same probability space
$(\Omega,\mathscr{F},P)$ with $\mathscr{F}=\sigma\{B\}\vee
\sigma\{W\}$. Let $(X_n)$ be a sequence of $\sigma\{B\}$-measurable
random variables, and let $X$ be a $\mathscr{F}$-measurable
random variable. In the sequel, we will write $X_n
\overset{{\rm stably}}{\longrightarrow}
X$ if $(Z,X_n)\overset{{\rm law}}{\longrightarrow}(Z,X)$ 
for all bounded and $\sigma\{B\}$-measurable
random variable $Z$.
In particular, we see that the stable convergence imply the convergence
in law. Moreover, it is easily checked that the convergence in probability 
implies the stable convergence.
We refer to \cite{JS}
for an exhaustive study of this notion.

Now, let us introduce the following object:
\begin{defi}\label{def1}
Let $f:\R^2\to\R$ be a continuously differentiable function, and fix 
a time $t>0$. 
Provided it exists, we define $\int_0^t \nabla f(B_s)\cdot dB_s$
to be the limit in probability, as $n\to\infty$, of
\begin{eqnarray}
I_n(t)&:=&
\sum_{k=0}^{\lfloor nt\rfloor -1} 
\frac{\frac{\partial f}{\partial x}(B^{(1)}_{k/n},B^{(2)}_{k/n}) 
+\frac{\partial f}{\partial x}(B^{(1)}_{(k+1)/n},B^{(2)}_{k/n}) 
}2
\big(B^{(1)}_{(k+1)/n}-B^{(1)}_{k/n})\notag\\
&&+\sum_{k=0}^{\lfloor nt\rfloor -1} \frac{\frac{\partial f}{\partial y}(B^{(1)}_{k/n},B^{(2)}_{k/n}) 
+\frac{\partial f}{\partial y}(B^{(1)}_{k/n},B^{(2)}_{(k+1)/n}) 
}2
\big(B^{(2)}_{(k+1)/n}-B^{(2)}_{k/n}).
\label{formula}
\end{eqnarray}
If $I_n(t)$ defined by (\ref{formula}) does not converge in probability but converges
stably, we denote the limit by
$\int_0^t \nabla f(B_s)\cdot d^\star B_s$. 
\end{defi}
Our main result is as follows:
\begin{thm}\label{main-thm}
Let $f:\R^2\to\R$ be a function belonging to $\mathscr{C}^8$ and verifying $({\bf H}_{\bf 8})$, see (\ref{hyp}) below.
Let also $B=(B^{(1)},B^{(2)})$
denote a 2D fractional Brownian motion of Hurst index $H\in(0,1)$, and $t>0$ be a fixed time.
\begin{enumerate}
\item If $H>1/4$ then $\int_0^t \nabla f(B_s)\cdot dB_s$ is well-defined, and we have
\begin{equation}\label{h>1/4}
f(B_t) = f(0) +\int_0^t \nabla f(B_s)\cdot dB_s. 
\end{equation}
\item If $H=1/4$ then only 
$\int_0^t \nabla f(B_s)\cdot d^{\star}B_s$ is well-defined, and we have
\begin{equation}\label{change}
f(B_t) \overset{{\rm Law}}{=} f(0) +
\int_0^t \nabla f(B_s)\cdot d^\star B_s
 + \frac{\sigma_{1/4}}{\sqrt{2}}
\int_0^t \frac{\partial^2 f}{\partial x\partial y}(B_s)dW_s. 
\end{equation}
Here, $\sigma_{1/4}$ is the universal constant defined below 
by (\ref{sigma}), and
$\int_0^t \nicefrac{\partial^2 f}{\partial x\partial y}(B_s)dW_s$ denotes a classical Wiener integral
with respect to the independent Brownian motion $W$.
\item If $H<1/4$ then the integral $\int_0^t B_s\cdot d^\star B_s$ does not exist. Therefore,
it is not possible to write a change of variable formula for $B^{(1)}_tB^{(2)}_t$
using the integral defined in Definition \ref{def1}. 
\end{enumerate}
\end{thm}

\begin{rem}\label{comments}
{\rm
\begin{enumerate}
\item Due to the definition of the stable convergence,
we can freely move each component in (\ref{change}) from the right hand side 
to the left (or
from the left hand side to the right). 
\item Whenever $\beta$ denotes a one-dimensional fractional Brownian motion
with Hurst index in $(0,1/2)$, it is easily checked, for any fixed $t>0$, that 
$\sum_{k=0}^{\lfloor nt\rfloor-1} \beta_{k/n}\big(\beta_{(k+1)/n}-\beta_{k/n}\big)$
does not converge in law. 
$\big($Indeed, on one hand, we have 
$$
\beta_{\lfloor nt\rfloor/t}^2
= \sum_{k=0}^{\lfloor nt\rfloor-1} \big(\beta_{(k+1)/n}^2 - \beta_{k/n}^2\big)= 
2\sum_{k=0}^{\lfloor nt\rfloor-1}\beta_{k/n}\big(\beta_{(k+1)/n}-\beta_{k/n}\big)
+\sum_{k=0}^{\lfloor nt\rfloor-1}\big(\beta_{(k+1)/n}-\beta_{k/n}\big)^2
$$
and, on the other hand, it is well-known (see e.g. \cite{KG}) that
$$
n^{2H-1}\sum_{k=0}^{\lfloor nt\rfloor-1} \big(\beta_{(k+1)/n}-\beta_{k/n}\big)^2 \underset{n\to\infty}{\overset{L^2}{\longrightarrow}} t.
$$
These two facts imply immediately that
$$
\sum_{k=0}^{\lfloor nt\rfloor-1}\beta_{k/n}\big(\beta_{(k+1)/n}-\beta_{k/n}\big)=
\frac12\left( \beta_{\lfloor nt\rfloor/t}^2 -
\sum_{k=0}^{\lfloor nt\rfloor-1}\big(\beta_{(k+1)/n}-\beta_{k/n}\big)^2\right)
$$
does not converge in law$\big)$. 
On the other hand, whenever $H>1/6$, the quantity 
$$\sum_{k=0}^{\lfloor nt\rfloor-1} \frac12\big(f(\beta_{k/n})+f(\beta_{(k+1)/n})\big)\,\big(\beta_{(k+1)/n}-\beta_{k/n}\big)$$
converges in $L^2$ for any regular enough function $f:\R\to\R$, see \cite{GNRV} and \cite{CN}.
This last fact roughly explains why there is a ``symmetric'' part
in the Riemann sum (\ref{formula}). 
\item We stress that it is still an open problem to know if each individual integral 
$\int_0^t \frac{\partial f}{\partial x}
(B_s)d^{(\star)}B^{\small(1\small)}_s$ and  $\int_0^t \frac{\partial f}{\partial y}
(B_s)d^{\small(\star\small)}B^{(2)}_s$ could be defined separately. Indeed, in the first two points of
Theorem \ref{main-thm}, we
``only'' prove that their sum, that is $\int_0^t \nabla f(B_s)\cdot d^{\small(\star\small)}B_s$, is well-defined.
\item Let us give a quicker proof of (\ref{change})
in the particular case where $f(x,y)=xy$. Let $\beta$ be a one-dimensional fractional Brownian motion of index $1/4$.
The classical Breuer-Major's theorem \cite{breuer-major} 
yields:
\begin{equation}\label{BM2}
\frac{1}{\sqrt{n}} \sum_{k=0}^{\lfloor n\cdot\rfloor-1} \big(\sqrt{n}(\beta_{(k+1)/n}-\beta_{k/n})^2 - 1\big)
\overset{\rm Law}{=}
\frac{1}{\sqrt{n}} \sum_{k=0}^{\lfloor n\cdot\rfloor-1} \big((\beta_{k+1}-\beta_{k})^2 - 1\big)
\overset{\rm stably}{\underset{n\to\infty}{\longrightarrow}}
\sigma_{1/4}\, W.
\end{equation}
Here, the convergence is stable and holds in the 
Skorohod space $\mathscr{D}$ of c\`adl\`ag functions on $[0,\infty)$. Moreover,
$W$ still denotes a standard Brownian motion independent of $\beta$
(the independence is a
consequence of the central limit theorem for multiple stochastic integrals
proved in \cite{peccati-tudor})
 and 
the constant
$\sigma_{1/4}$ is given by 
\begin{equation}\label{sigma}
\sigma_{1/4}:=\sqrt{\frac{1}{2}\sum_{k\in\mathbb{Z}} \left(\sqrt{|k+1|}+\sqrt{|k-1|} - 2\sqrt{|k|}\right)^2}<\infty.
\end{equation}
Now, let $\widetilde{\beta}$ be another fractional Brownian motion of index $1/4$, independent of $\beta$. 
From (\ref{BM2}), we get
\begin{eqnarray*}
\left(
\frac{1}{\sqrt{n}} \sum_{k=0}^{\lfloor nt\rfloor-1} \big(\sqrt{n}(\beta_{(k+1)/n}-\beta_{k/n})^2 - 1\big),
\frac{1}{\sqrt{n}} \sum_{k=0}^{\lfloor nt\rfloor -1} \big(\sqrt{n}(\widetilde{\beta}_{(k+1)/n}-\widetilde{\beta}_{k/n})^2 - 1\big)
\right)
\\
\overset{\rm stably}{\underset{n\to\infty}{\longrightarrow}}
\sigma_{1/4}\,(W,\widetilde{W})
\end{eqnarray*}
for $(W,\widetilde{W})$ a 2D standard Brownian motion, independent of the 2D fractional Brownian motion
$(\beta,\widetilde{\beta})$. 
In particular, by difference, we have
$$
\frac12\sum_{k=0}^{\lfloor n\cdot\rfloor -1} \big((\beta_{(k+1)/n}-\beta_{k/n})^2 - 
(\widetilde{\beta}_{(k+1)/n}-\widetilde{\beta}_{k/n})^2\big)
\overset{\rm stably}{\underset{n\to\infty}{\longrightarrow}}
\frac{\sigma_{1/4}}{2}(W-\widetilde{W})
\overset{\rm Law}{=}
\frac{\sigma_{1/4}}{\sqrt{2}} \,W.
$$
Now, set $B^{(1)}=(\beta+\widetilde{\beta})/\sqrt{2}$ and $B^{(2)}=(\beta-\widetilde{\beta})/\sqrt{2}$. 
It is easily checked that $B^{(1)}$ and $B^{(2)}$ are two independent fractional Brownian
motions of index $1/4$. Moreover, we can rewrite the previous convergence as
\begin{equation}\label{BM5}
\sum_{k=0}^{\lfloor n\cdot\rfloor-1} (B^{(1)}_{(k+1)/n}-B^{(1)}_{k/n}) 
(B^{(2)}_{(k+1)/n}-B^{(2)}_{k/n})
\overset{\rm stably}{\underset{n\to\infty}{\longrightarrow}}
\frac{\sigma_{1/4}}{\sqrt{2}} \,W,
\end{equation}
with $B^{(1)}$, $B^{(2)}$ and $W$ independent.
On the other hand, for any $a,b,c,d\in\R$:
$$
bd - ac = a(d-c) + c(b-a) + (b-a)(d-c).
$$
Choosing $a=B^{(1)}_{k/n}$, $b=B^{(1)}_{(k+1)/n}$, $c=B^{(2)}_{k/n}$ and $d=B^{(2)}_{(k+1)/n}$, 
and suming for $k$ over $0,\ldots,\lfloor nt\rfloor -1$, we obtain
\begin{eqnarray}
B^{(1)}_{\lfloor nt\rfloor/n}\,B^{(2)}_{\lfloor nt\rfloor/n} &=& \sum_{k=0}^{\lfloor nt\rfloor -1} 
B^{(1)}_{k/n}\big(B^{(2)}_{(k+1)/n}-B^{(2)}_{k/n}\big)
+B^{(2)}_{k/n}\big(B^{(1)}_{(k+1)/n}-B^{(1)}_{k/n}\big)\notag\\
&&
+ \sum_{k=0}^{\lfloor nt\rfloor -1} 
\big(B^{(1)}_{(k+1)/n}-B^{(1)}_{k/n}\big)\big(B^{(2)}_{(k+1)/n}-B^{(2)}_{k/n}\big).\label{eheh}
\end{eqnarray}
Hence, passing to the limit using (\ref{BM5}), we get the desired conclusion 
in (\ref{change}), in the particular case where $f(x,y)=xy$.
Note that the second term in the right-hand side of (\ref{eheh}) is
the discrete analogue of the 2-{\it covariation}
introduced by Errami and Russo in \cite{errami-russo}.
\item We could prove (\ref{change})
at a functional level (note that it has precisely been done for $f(x,y)=xy$
in the proof just below). But, in order to keep the length of this paper within limits, 
we defer to future analysis this rather technical investigation.
\item In the very recent work \cite{nourdin-reveillac}, R\'eveillac and I  proved the following result (see also
Burdzy and Swanson \cite{burdzy-swanson} for similar results in the case where
$\beta$ is replaced by the solution of the
stochastic heat equation driven by a space/time white noise). If 
$\beta$ denotes a one-dimensional fractional Brownian motion of index $1/4$ and if $g:\R\to\R$
is regular enough, then
\begin{equation}\label{nr}
\frac{1}{\sqrt{n}}\sum_{k=0}^{n-1}
g(\beta_{k/n})\big(\sqrt{n}(\beta_{(k+1)/n}-\beta_{k/n})^2 -1\big)
\underset{n\to\infty}{\overset{\rm stably}{\longrightarrow}}
\frac14\int_0^1 g''(\beta_s)ds + \sigma_{1/4}\int_0^1 g(\beta_s)dW_s
\end{equation}
for $W$ a standard Brownian motion independent of $\beta$. Compare with 
Proposition \ref{law} below. In particular, by choosing $g$ identically one in (\ref{nr}), it agrees with (\ref{BM2}). 
\item The fractional Brownian motion of index $1/4$ has a remarkable physical interpretation
in terms of particle systems. Indeed, if one consider an infinite number of particles,
initially placed on the real line according to a Poisson distribution, performing
independent Brownian motions and undergoing ``elastic'' collisions, then the trajectory
of a fixed particle (after rescaling) converges to a fractional Brownian motion
of index $1/4$. See Harris \cite{harris} for heuristic arguments, and D\"urr, Goldstein and Lebowitz \cite{DGL}
for precise results.
\end{enumerate}
}
\end{rem}

Now, the rest of the note is entirely devoted to the proof of Theorem \ref{main-thm}.
The Section 2 contains some preliminaries and fix the notation.
Some technical results are postponed in Section 3. 
Finally, the proof of Theorem \ref{main-thm} is done in Section 4.

\section{Preliminaries and notation}\label{prelim}

We shall now provide a short description of the tools of Malliavin calculus that will be needed in the
following sections. The reader is referred to the monographs \cite{malliavin} and \cite{nualart} for
any unexplained notion or result.

Let $B=(B^{(1)}_t,B^{(2)}_t)_{t\in[0,T]}$ be a 2D fractional Brownian motion with Hurst parameter
belonging to $(0,1/2)$. We denote by $\mathcal{H}$ the Hilbert space defined as the closure of the
set of step $\mathbb{R}^2$-valued functions on $[0,T]$, with respect to the scalar product induced
by
$$
\left\langle \big({\bf 1}_{[0,t_1]},{\bf 1}_{[0,t_2]}\big),\big({\bf 1}_{[0,s_1]},{\bf 1}_{[0,s_2]}\big)\right\rangle
_\mathcal{H} = R_H(t_1,s_1)+ R_H(t_2,s_2),\quad s_i,t_i\in[0,T],\quad i=1,2,
$$
where $R_H(t,s)=\frac{1}{2}\left( t^{2H}+s^{2H}-|t-s|^{2H}\right)$. 
The mapping $({\bf 1}_{[0,t_1]},{\bf 1}_{[0,t_2]})\mapsto B_{t_1}^{(1)}+B_{t_2}^{(2)}$
can be extended to an isometry between $\mathcal{H}$ and the Gaussian space
associated with $B$.
Also, $\mathfrak{H}$ will denote
the Hilbert space defined as the closure of the
set of step $\mathbb{R}$-valued functions on $[0,T]$, with respect to the scalar product induced
by
$$
\left\langle {\bf 1}_{[0,t]},{\bf 1}_{[0,s]}\right\rangle
_\mathfrak{H} = R_H(t,s),\quad s,t\in[0,T].
$$
The mapping ${\bf 1}_{[0,t]}\mapsto B_{t}^{(i)}$ ($i$ equals $1$ or $2$)
can be extended to an isometry between $\mathfrak{H}$ and the Gaussian space
associated with $B^{(i)}$.

Consider the set of all smooth cylindrical random variables, i.e. of the form
\begin{equation}\label{eq:cylindrical}
F=f\big(B(\varphi_1),\ldots,B(\varphi_k)\big),\quad\varphi_i\in\mathcal{H},\quad i=1,\ldots,k,
\end{equation}
where $f\in\mathscr{C}^{\infty}$ is bounded with bounded
derivatives. The derivative operator $D$ of a smooth cylindrical random variable of
the above form is defined as the $\mathcal{H}$-valued random variable
$$
DF = \sum_{i=1}^{k}\frac{\partial f}{\partial x_i}\big(B(\varphi_1),\ldots,B(\varphi_k)\big)\varphi_i=:\big( D_{B^{(1)}}F,D_{B^{(2)}}F \big).
$$
In particular, we have
$$
D_{B^{(i)}}B^{(j)}_t = \delta_{ij}{\bf 1}_{[0,t]}\quad\mbox{for $i,j\in\{1,2\}$, and $\delta_{ij}$ the
Kronecker symbol.}
$$
By iteration, one can define the $m$th derivative $D^m F$ (which is a
symmetric element of $L^2(\Omega,\mathcal{H}^{\otimes m}))$
for $m\geq 2$. As usual, for any $m\geq 1$, the space $\mathbb{D}^{m,2}$ denotes the closure of the set
of smooth random variables with respect to the norm $\|\cdot\|_{m,2}$ defined by the relation
$$
\| F\|^2_{m,2} = E| F|^2 + \sum_{i=1}^m E\|D^i F\|^2_{\mathcal{H}^{\otimes i}}.
$$
The derivative $D$ verifies the chain rule. Precisely, if $\varphi:\R^n\to\R$ belongs to $\mathscr{C}^1$ with bounded derivatives
and if $F_i$, $i=1,\ldots,n$, are in $\mathbb{D}^{1,2}$, then $\varphi(F_1,\ldots,F_n)\in\mathbb{D}^{1,2}$ and
$$D\varphi(F_1,\ldots,F_n)=\sum_{i=1}^{n} \frac{\partial \varphi}{\partial x_i}(F_1,\ldots,F_n)DF_i.$$
The $m$th derivative $D_{B^{(i)}}^m$ ($i$ equals $1$ or $2$) verifies the following Leibnitz rule: for any $F,G\in\mathbb{D}^{m,2}$ such that
$FG\in\mathbb{D}^{m,2}$, we have
\begin{equation}\label{leibnitz}
\big(D_{B^{(i)}}^{m}FG\big)_{t_1,\ldots,t_m}= \sum 
\big(D_{B^{(i)}}^{r}F\big)_{s_1,\ldots,s_r}\big(D_{B^{(i)}}^{m-r}G\big)_{u_1,\ldots,u_{m-r}},\quad t_i\in[0,T],\quad i=1,\ldots,m,
\end{equation}
where the sum runs over any subset $\{s_1,\ldots,s_r\}\subset\{t_1,\ldots,t_m\}$
and where we write $\{ t_1,\ldots,t_m\}\setminus\{s_1,\ldots,s_r\}=:
\{u_1,\ldots,u_{m-r}\}$.

The divergence operator $\delta$ is the adjoint of the derivative operator. If a random variable
$u\in L^{2}(\Omega,\mathcal{H})$ belongs to ${\rm dom}\delta$, the domain of the divergence operator,
then $\delta(u)$ is defined by the duality relationship
$$
E\big(F\delta(u)\big)=E\langle DF,u\rangle_\mathcal{H}
$$
for every $F\in\mathbb{D}^{1,2}$. 

For every $q\geq 1$, let $\mathcal{H}_{q}$ be the $q$th Wiener chaos of $B$,
that is, the closed linear subspace of $L^{2}\left( \Omega ,\mathcal{A}
,P\right) $ generated by the random variables $\{H_{q}\left( B\left(
h\right) \right) ,h\in \mathcal{H},\| h\| _{\mathcal{H}}=1\}$, where $H_{q}$
is the $q$th Hermite polynomial given by $H_q(x)=(-1)^q e^{x^2/2} \frac{d^q}{dx^q}\big(e^{-x^2/2}\big)$. The mapping 
\begin{equation}\label{hermite}
I_{q}(h^{\otimes
q})=H_{q}\left( B\left( h\right) \right)
\end{equation} 
provides a linear isometry
between the symmetric tensor product $\mathcal{H}^{\odot q}$ 
(equipped with the modified norm $\frac1{\sqrt{q!}}\|\cdot\|_{\mathcal{H}^{\otimes q}}$)
and $\mathcal{H}
_{q}$. 
The following duality formula holds%
\begin{equation}\label{banff}
E\left( FI_{q}(f)\right) =E\left( \left\langle D^{q}F,f\right\rangle _{
\mathcal{H}^{\otimes q}}\right) ,
\end{equation}
for any $f\in \mathcal{H}^{\odot q}$ and  $F\in 
\mathbb{D}^{q,2}$. In particular, we have
\begin{equation}
E\left( F I^{(i)}_{q}(g)\right) =E\left( \left\langle D_{B^{(i)}}^{q}F,g\right\rangle _{
\mathfrak{H}^{\otimes q}}\right) ,\quad i=1,2,  \label{dual2}
\end{equation}
for any $g\in \EuFrak H^{\odot q}$ and  $F\in 
\mathbb{D}^{q,2}$, where, for simplicity, we write $I_q^{(i)}(g)$ 
whenever  the corresponding $q$th multiple integral is only with respect to
$B^{(i)}$. 

Finally, we mention the following particular case (actually, the only one we will need in the sequel) of the classical
multiplication formula: if $f,g\in \HH$, $q\geq 1$ and $i\in\{1,2\}$,
then
\begin{equation}\label{multiplication}
I^{(i)}_q(f^{\otimes q})I^{(i)}_q(g^{\otimes q})=\sum_{r=0}^{q} r!\binom{q}{r}^2
I^{(i)}_{2q-2r}(f^{\otimes q-r}\otimes g^{\otimes q-r})\langle f,g\rangle_\HH^r.
\end{equation}

\section{Some technical results}

In this section, we collect some crucial results for the proof of (\ref{change}),
the only case which is difficult.

Here and in the rest of the paper, we set 
$$\Delta B^{(i)}_{k/n} := B^{(i)}_{(k+1)/n} - B^{(i)}_{k/n},\quad 
\delta_{k/n}:={\bf 1}_{[k/n,(k+1)/n]}\quad\mbox{and}\quad
\e_{k/n}:={\bf 1}_{[0,k/n]},$$
for any $i\in\{1,2\}$ and $k\in\{0,\ldots,n-1\}$.

In the sequel, for $g:\R^2\to\R$ belonging to $\mathscr{C}^q$, we will need assumption of the type:\\
\begin{equation}\label{hyp}
({\bf H}_{\bf q})\quad
\sup_{s\in[0,1]} E\left|\frac{\partial^{a+b}g}{\partial x^a\partial y^b}(B^{(1)}_s,B^{(2)}_s) \right|^p<\infty\,
\mbox{ for all $p\geq 1$ and all integers $a,b\geq 0$ s.t. $a+b\leq q$}.
\end{equation}

We begin by the following technical lemma:
\begin{lemme}
\label{lemma:Tec1}
Let $\beta$ be a 1D fractional Brownian motion of Hurst index $1/4$. We have
\begin{itemize}
\item[(i)] $\left\vert E\big(\beta_r (\beta_t-\beta_s)\big)\right\vert \leq \sqrt{|t-s|}$ for any 
$0\leq r,s,t\leq 1$,
\item[(ii)] $\displaystyle{\sum_{k,l=0}^{n-1} \left\vert \left\langle \varepsilon_{l/n}, 
\delta_{k/n} \right\rangle_{\EuFrak H} \right\vert \underset{n\to\infty}{=} O(n)},$
\item[(iii)] $\displaystyle{
\sum_{k,l=0}^{n-1} \left\vert \left\langle \delta_{l/n}, 
\delta_{k/n} \right\rangle_{\EuFrak H} \right\vert^r
 \underset{n\to\infty}{=} O(n^{1-r/2})}$ for any $r\geq 1$,
\item[(iv)] $\displaystyle{\sum_{k=0}^{n-1}\left\vert \left\langle \varepsilon_{k/n}, 
\delta_{k/n} \right\rangle_{\EuFrak H} + \frac{1}{2 \sqrt{n}} \right\vert  
\underset{n\to\infty}{=} O(1)}$,
\item[(v)] $\displaystyle{\sum_{k=0}^{n-1}\left\vert \left\langle \varepsilon_{k/n}, 
\delta_{k/n} \right\rangle_{\EuFrak H}^2 - \frac{1}{4 n} \right\vert  
\underset{n\to\infty}{=} O(1/\sqrt{n})}$.
\end{itemize}
\end{lemme}
\noindent
{\bf Proof of Lemma \ref{lemma:Tec1}}. 
\begin{itemize}
\item[\emph{(i)}] We have%
$$
E\big(\beta_{r}(\beta_{t}-\beta_{s})\big) 
=\frac{1}{2}\big(\sqrt{t}-\sqrt{s}\big)+\frac{1}{2}\left( \sqrt{|s-r|}-\sqrt{|t-r|}\right).
$$
Using the classical inequality $\big|\sqrt{|b|}-\sqrt{|a|}\big|\le
\sqrt{|b-a|}$, the desired result follows. 
\item[\emph{(ii)}] Observe that 
$$ \left\langle \varepsilon_{l/n},\delta_{k/n} \right\rangle_{\EuFrak H}= 
\frac{1}{2 \sqrt{n}} \left( \sqrt{k+1}-\sqrt{k}-\sqrt{\vert k+1-l\vert} + 
\sqrt{\vert k- l \vert} \right).$$
Consequently, for any fixed $l\in\{0,\ldots,n-1\}$, we have
\begin{eqnarray*}
\sum_{k=0}^{n-1}
\left|\left\langle \varepsilon_{l/n},\delta_{k/n} \right\rangle_{\EuFrak H}\right|
&\le&\frac12+\frac1{2\sqrt{n}}\left(
\sum_{k=0}^{l -1}
\sqrt{l-k}-\sqrt{l-k-1}
\right.\\
&&\left.
+1+
\sum_{k=l +1}^{n -1}
\sqrt{k-l+1}-\sqrt{k-l}
\right)\\
&=&\frac12 + \frac{1}{2\sqrt{n}}\big(\sqrt{l}+\sqrt{n-l}\big)
\end{eqnarray*}
from which we deduce that
$\displaystyle{\sup_{0\leq l\leq n-1} \sum_{k=0}^{n-1} \left\vert \left\langle \varepsilon_{l/n}, 
\delta_{k/n} \right\rangle_{\EuFrak H} \right\vert \underset{n\to\infty}{=}O(1)}$.
It follows that
$$
\sum_{k,l=0}^{n-1}  \left\vert \left\langle \varepsilon_{l/n}, \delta_{k/n} \right\rangle_{\EuFrak H} 
\right\vert \leq n \sup_{0\leq l\leq n-1} \sum_{k=0}^{n-1} \left\vert \left\langle \varepsilon_{l/n}, 
\delta_{k/n} \right\rangle_{\EuFrak H} \right\vert 
\underset{n\to\infty}{=} O(n).
$$
\item[\emph{(iii)}]
We have, by noting $\rho(x)=\frac12\big(\sqrt{|x+1|}+\sqrt{|x-1|}-2\sqrt{|x|}\big)$:
\begin{eqnarray*}
\sum_{k,l=0}^{n-1} \left\vert \left\langle \delta_{l/n}, 
\delta_{k/n} \right\rangle_{\EuFrak H} \right\vert^r &=&
n^{-r/2}\sum_{k,l=0}^{n-1} \left\vert \rho^r(l-k) \right\vert
\leq n^{1-r/2}\sum_{k\in\Z}\left\vert \rho^r(k) \right\vert.
\end{eqnarray*} 
Since $\sum_{k\in\Z}\left\vert \rho^r(k) \right\vert<\infty$ if $r\geq 1$, the desired conclusion follows.
\item[\emph{(iv)}] 
is a consequence of the following identity combined with a telescopic sum argument:
$$\left\vert \left\langle \varepsilon_{k/n}, \delta_{k/n} \right\rangle_{\EuFrak H} + \frac{1}{2\sqrt{n}} 
\right\vert=\frac{1}{2\sqrt{n}}\big(\sqrt{k+1}-\sqrt{k}\big).
$$
\item[\emph{(v)}] 
We have
$$\left\vert \left\langle \varepsilon_{k/n}, \delta_{k/n} \right\rangle_{\EuFrak H}^2 - \frac{1}{4 n} 
\right\vert=\frac{1}{4 n} \left(\sqrt{k+1}-\sqrt{k}\right) \left| \sqrt{k+1}-\sqrt{k}-2 \right|.$$
Thus, the desired bound is immediately checked 
by combining a telescoping sum argument with the fact that
$$
\left| \sqrt{k+1}-\sqrt{k}-2 \right| = \left| \frac1{\sqrt{k+1}+\sqrt{k}} - 2\right|\leq 2.
$$
\end{itemize}
\fin

Also the following lemma will be useful in the sequel:
\begin{lemme}\label{lm1}
Let $\alpha\geq 0$ and $q\geq 2$ be two positive integers, $g:\R^2\to\R$ be any function
belonging to $\mathscr{C}^{2q}$ and verifying $({\bf H}_{\bf 2q})$ defined by (\ref{hyp}),
and $B=(B^{(1)},B^{(2)})$ be a 2D fractional Brownian motion of Hurst index $1/4$.
Set
$$
V_n=n^{-q/4}\sum_{k=0}^{n-1} g(B^{(1)}_{k/n},B^{(2)}_{k/n})
\big(\Delta B^{(1)}_{k/n}\big)^\alpha \,
H_q\big(n^{1/4}\Delta B^{(2)}_{k/n}\big),
$$
where $H_q$ denotes the $q$th Hermite polynomial defined by $H_q(x)=(-1)^q e^{x^2/2} \frac{d^q}{dx^q}\big(e^{-x^2/2}\big)$.
Then, the following bound is in order: 
\begin{equation}\label{bound}
E\left(|V_n|^2\right) = O(n^{1-q/2-\alpha/2})\quad\mbox{ as }n\to\infty.
\end{equation}
\end{lemme}
{\bf Proof of Lemma \ref{lm1}}. 
We can write
\begin{eqnarray}
E\left(| V_n|^2\right) &=& n^{-q/2}\sum_{k,l=0}^{n-1}
E\big[
g(B^{(1)}_{k/n},B^{(2)}_{k/n})g(B^{(1)}_{l/n},B^{(2)}_{l/n})
\big(\Delta B^{(1)}_{k/n}\big)^\alpha \,\big(\Delta B^{(1)}_{l/n}\big)^\alpha \notag\\
&&\hskip6cm \times
H_q\big(n^{1/4}\Delta B^{(2)}_{k/n}\big)H_q\big(n^{1/4}\Delta B^{(2)}_{l/n}\big)
\big]\notag\\
&\underset{(\ref{hermite})}{=}& \sum_{k,l=0}^{n-1}
E\big[
g(B^{(1)}_{k/n},B^{(2)}_{k/n})g(B^{(1)}_{l/n},B^{(2)}_{l/n})
\big(\Delta B^{(1)}_{k/n}\big)^\alpha \,\big(\Delta B^{(1)}_{l/n}\big)^\alpha 
I^{(2)}_q(\delta_{k/n}^{\otimes q}) I^{(2)}_q(\delta_{l/n}^{\otimes q})\big]\notag\\
&\underset{(\ref{multiplication})}{=}& \sum_{r=0}^q r!\binom{q}{r}^2 \,\,\sum_{k,l=0}^{n-1}
E\big[
g(B^{(1)}_{k/n},B^{(2)}_{k/n})g(B^{(1)}_{l/n},B^{(2)}_{l/n})\notag\\
&&\hskip1cm\times 
\big(\Delta B^{(1)}_{k/n}\big)^\alpha \,\big(\Delta B^{(1)}_{l/n}\big)^\alpha 
I^{(2)}_{2q-2r}(\delta_{k/n}^{\otimes q-r}\otimes \delta_{l/n}^{\otimes q-r})\big]\langle 
\delta_{k/n},\delta_{l/n}\rangle_\HH^r\notag\\
&\underset{(\ref{dual2})}{=}& \sum_{r=0}^q r!\binom{q}{r}^2 \sum_{k,l=0}^{n-1}
E\left\langle
D_{B^{(2)}}^{2q-2r}\left(g(B^{(1)}_{k/n},B^{(2)}_{k/n})g(B^{(1)}_{l/n},B^{(2)}_{l/n})\right.\right.\notag\\
&&\hskip1cm\times\left. \left.
\big(\Delta B^{(1)}_{k/n}\big)^\alpha \,\big(\Delta B^{(1)}_{l/n}\big)^\alpha\right), 
\delta_{k/n}^{\otimes q-r}\otimes \delta_{l/n}^{\otimes q-r}\right\rangle_{\HH^{\otimes 2q-2r}}
\langle\delta_{k/n},\delta_{l/n}\rangle_\HH^r\notag\\
&\underset{(\ref{leibnitz})}{=}& \sum_{r=0}^q r!\binom{q}{r}^2 \sum_{a+b=2q-2r} \frac{(a+b)!}{a!b!}\sum_{k,l=0}^{n-1}
E\left(\frac{d^a g}{dy^a}(B^{(1)}_{k/n},B^{(2)}_{k/n})\frac{d^b g}{dy^b}(B^{(1)}_{l/n},B^{(2)}_{l/n})\right.\notag\\
&&\times\left.
\big(\Delta B^{(1)}_{k/n}\big)^\alpha \,\big(\Delta B^{(1)}_{l/n}\big)^\alpha\right)
(2q-2r)!\left\langle \e_{k/n}^{\otimes a}\widetilde{\otimes}\e_{l/n}^{\otimes b}, 
\delta_{k/n}^{\otimes q-r}\otimes \delta_{l/n}^{\otimes q-r}\right\rangle_{\HH^{\otimes 2q-2r}}
\langle\delta_{k/n},\delta_{l/n}\rangle_\HH^r.\notag\\
\label{cracboumue}
\end{eqnarray}
Now, observe that, uniformly in $k,l\in\{0,\ldots,n-1\}$:
$$
\left\langle \e_{k/n}^{\otimes a}\widetilde{\otimes}\e_{l/n}^{\otimes b}, 
\delta_{k/n}^{\otimes q-r}\otimes \delta_{l/n}^{\otimes q-r}\right\rangle_{\HH^{\otimes 2q-2r}}
\underset{n\to\infty}{=}
O(n^{-(q-r)}),\,\mbox{ see Lemma \ref{lemma:Tec1} (i)},
$$
$$
\left|
E\left(\frac{d^a g}{dy^a}(B^{(1)}_{k/n},B^{(2)}_{k/n})\frac{d^b g}{dy^b}(B^{(1)}_{l/n},B^{(2)}_{l/n})
\big(\Delta B^{(1)}_{k/n}\big)^\alpha \,\big(\Delta B^{(1)}_{l/n}\big)^\alpha\right)\right|
\underset{n\to\infty}{=}O(n^{-\alpha/2}),\,\mbox{use $({\bf H}_{\bf 2q})$},
$$
and, also: 
$$
\sum_{k,l=0}^{n-1} \langle \delta_{k/n},\delta_{l/n}\rangle_\HH^r = O(n^{1-r/2})\quad\mbox{for any fixed $r\geq 1$, see Lemma \ref{lemma:Tec1} (iii)}.
$$
Finally, the desired conclusion is obtained by plugging these three bounds 
into (\ref{cracboumue}), after having separated the cases $r=0$ and $r=1$.
\fin 

The independent Brownian motion appearing in (\ref{change}) comes from the following proposition. 

\begin{prop}\label{law}
Let $(\beta,\widetilde{\beta})$ be a 2D fractional Brownian motion of Hurst index $1/4$.
Consider two functions $g,\widetilde{g}:\R^2\to\R$ belonging in $\mathscr{C}^4$, and  assume that they both verify $({\bf H}_{\bf 4})$
defined by (\ref{hyp}). Then
\begin{eqnarray*}
&&(G_n,\widetilde{G}_n):=\left(
\frac{1}{\sqrt{n}}
\sum_{k=0}^{n-1} g(\beta_{k/n},\widetilde{\beta}_{k/n})\big(\sqrt{n}(\Delta \beta_{k/n})^2 - 1\big),
\frac{1}{\sqrt{n}}
\sum_{k=0}^{n-1} \widetilde{g}(\beta_{k/n},\widetilde{\beta}_{k/n})\big(\sqrt{n}(\Delta \widetilde{\beta}_{k/n})^2 - 1\big)
\right)\\
&&\underset{n\to\infty}{\overset{\rm stably}{\longrightarrow}}
\left(
\sigma_{1/4}\int_0^1 g(\beta_s,\widetilde{\beta}_s)dW_s
+\frac14\int_0^1 \frac{\partial^2 g}{\partial x^2}(\beta_s,\widetilde{\beta}_s)ds,
\sigma_{1/4}\int_0^1 \widetilde{g}(\beta_s,\widetilde{\beta}_s)d\widetilde{W}_s
+\frac14
\int_0^1 \frac{\partial^2 \widetilde{g}}{\partial y^2}(\beta_s,\widetilde{\beta}_s)ds\right),
\end{eqnarray*}
where $(W,\widetilde{W})$ is a 2D standard Brownian motion independent of $(\beta,\widetilde{\beta})$, and
$\sigma_{1/4}$ is defined by (\ref{sigma}).
\end{prop}
In the particular case where $g(x,y)=g(x)$ and $
\widetilde{g}
(x,y)=\widetilde{g}(y)$, the conclusion of the proposition
follows directly from (\ref{nr}). In the general case, the proof
only consists to extend literaly the proof of (\ref{nr}) contained
in \cite{nourdin-reveillac}.
Details are left to the reader.

\section{Proof of Theorem \ref{main-thm}}
We are now in position to prove our main result, that is
Theorem \ref{main-thm}.\\
\\
{\bf Proof of the third point} (case $H<1/4$). 
Firstly, observe that (\ref{BM2}) is actually a particular case of the following result, which
is valid for any fractional Brownian $\beta$ with Hurst index $H$
belonging to $(0,3/4)$:
$$
\frac{1}{\sqrt{n}} \sum_{k=0}^{\lfloor n\cdot\rfloor-1} \big(
n^{2H}(\Delta\beta_{k/n})^2 - 1\big)
\overset{\rm stably}{\underset{n\to\infty}{\longrightarrow}}
\sigma_{H}\, W
$$
with $W$ an independent Brownian motion and $\sigma_H>0$ an (explicit) constant.
By mimicking the proof contained in the fourth point of 
Remark 1.3, we get, here, for any $H\in(0,3/4)$,
\begin{equation}\label{vla}
n^{2H-1/2}\sum_{k=0}^{\lfloor n\cdot\rfloor-1} \Delta B^{(1)}_{k/n}\Delta B^{(2)}_{k/n} \overset{\rm stably}{\underset{n\to\infty}{\longrightarrow}} \frac{\sigma_H}{\sqrt{2}}W.
\end{equation}
But, see (\ref{eheh}), the existence of $\int_0^\cdot B_s\cdot d^\star B_s$ would imply in particular that
$\sum_{k=0}^{\lfloor n\cdot\rfloor-1} \Delta B^{(1)}_{k/n}\Delta B^{(2)}_{k/n}$
converges in law as $n\to\infty$, which is in contradiction with (\ref{vla}) for $H<1/4$.
The proof of the third point is done.\\
\\
{\bf Proof of the second point} (case $H=1/4$).
For the simplicity of the exposition, we assume from now 
that $t=1$, the general case
being of course similar up to cumbersome notation. 
For any $a,b,c,d\in\R$, by the classical Taylor formula,
we can expand $f(b,d)$ as (compare with (\ref{eheh})):
\begin{eqnarray}
&&f(a,c) + \partial_1 f(a,c)(b-a) +\partial_2 f(a,c) (d-c) + \frac12 \partial_{11}f(a,c)(b-a)^2
+\frac12\partial_{22} f(a,c)(d-c)^2\notag\\
&+&\frac1{6}\partial_{111} f(a,c) (b-a)^3 +\frac1{6}\partial_{222} f(a,c) (d-c)^3 
+ \frac1{24}\partial_{1111}f(a,c)(b-a)^4 + \frac1{24}\partial_{2222}f(a,c)(d-c)^4\notag\\
&+&\partial_{12}f(a,c)(b-a)(d-c)+\frac{1}{2}\partial_{112}f(a,c)(b-a)^2(d-c) + \frac12\partial_{122}f(a,c)(b-a)(d-c)^2\notag\\
&+&\frac16\partial_{1112}f(a,c)(b-a)^3(d-c)+\frac14\partial_{1122}f(a,c)(b-a)^2(d-c)^2+\frac16\partial_{1222}f(a,c)(b-a)(d-c)^3\notag\\
\label{ito2}
\end{eqnarray}
plus a remainder term. Here, as usual, the notation $\partial_{1\ldots 12\ldots 2}f$ (where the index 1 is repeated $k$
 times and the index 2 is repeated $l$ times) means that $f$ is differentiated $k$ times w.r.t. 
the first component and $l$ times w.r.t. the second one.
By combining (\ref{ito2}) with the following identity, available for any $h:\R\to\R$ belonging to $\mathscr{C}^4$:
\begin{eqnarray*}
&&h'(a)(b-a)+\frac12 h''(a)(b-a)^2 + \frac16 h'''(a)(b-a)^3 + \frac1{24} h''''(a)(b-a)^4\\
&=&\frac{h'(a)+h'(b)}{2}(b-a)-\frac1{12}h'''(a)(b-a)^3
-\frac1{24}h''''(a)(b-a)^4 + \mbox{ some remainder}
\end{eqnarray*}
we get that $f(b,d)$ can also be expanded
as
\begin{eqnarray}
&&f(a,c) + \frac12\big(\partial_1 f(a,c)+\partial_1 f(b,c)\big)(b-a) 
-\frac1{12}\partial_{111} f(a,c) (b-a)^3 - \frac1{24}\partial_{1111}f(a,c)(b-a)^4\notag\\
&+&\frac12\big(\partial_2 f(a,c)+\partial_2 f(a,d)\big)(d-c) 
-\frac1{12}\partial_{222} f(a,c) (d-c)^3 - \frac1{24}\partial_{2222}f(a,c)(d-c)^4\notag\\
&+&\partial_{12}f(a,c)(b-a)(d-c)+\frac{1}{2}\partial_{112}f(a,c)(b-a)^2(d-c) + \frac12\partial_{122}f(a,c)(b-a)(d-c)^2\notag\\
&+&\frac16\partial_{1112}f(a,c)(b-a)^3(d-c)+\frac14\partial_{1122}f(a,c)(b-a)^2(d-c)^2+\frac16\partial_{1222}f(a,c)(b-a)(d-c)^3\notag\\
\label{ito}
\end{eqnarray}
plus a remainder term. 

By setting $a=B^{(1)}_{k/n}$, $b=B^{(1)}_{(k+1)/n}$, $c=B^{(2)}_{k/n}$ and $d=B^{(2)}_{(k+1)/n}$
in (\ref{ito}), 
and by suming the obtained expression for $k$ over $0,\ldots,n -1$, we
deduce that the conclusion in Theorem \ref{main-thm} is a consequence of the following 
convergences:
\begin{eqnarray}
S_n^{(1)}&:=&\sum_{k=0}^{n-1} \partial_{111} f(B^{(1)}_{k/n},B^{(2)}_{k/n}) \big(\Delta B^{(1)}_{k/n}\big)^3
\underset{n\to\infty}{\overset{L^2}{\longrightarrow}}-\frac32\int_0^1 \partial_{1111} f(B^{(1)}_s,B^{(2)}_s)ds\label{sn1}\\
S_n^{(2)}&:=&\sum_{k=0}^{n-1} \partial_{1111} f(B^{(1)}_{k/n},B^{(2)}_{k/n}) \big(\Delta B^{(1)}_{k/n}\big)^4
\underset{n\to\infty}{\overset{L^2}{\longrightarrow}}3\int_0^1 \partial_{1111} f(B^{(1)}_s,B^{(2)}_s)ds\label{sn2}\\
S_n^{(3)}&:=&\sum_{k=0}^{n-1} \partial_{222} f(B^{(1)}_{k/n},B^{(2)}_{k/n}) \big(\Delta B^{(2)}_{k/n}\big)^3
\underset{n\to\infty}{\overset{L^2}{\longrightarrow}}-\frac32\int_0^1 \partial_{2222} f(B^{(1)}_s,B^{(2)}_s)ds\label{sn3}\\
S_n^{(4)}&:=&\sum_{k=0}^{n-1} \partial_{2222} f(B^{(1)}_{k/n},B^{(2)}_{k/n}) \big(\Delta B^{(2)}_{k/n}\big)^4
\underset{n\to\infty}{\overset{L^2}{\longrightarrow}}3\int_0^1 \partial_{2222} f(B^{(1)}_s,B^{(2)}_s)ds\label{sn4}\\
S_n^{(5)}&:=&\sum_{k=0}^{n-1} \partial_{12} f(B^{(1)}_{k/n},B^{(2)}_{k/n}) \Delta B^{(1)}_{k/n}\,\Delta B^{(2)}_{k/n}
\underset{n\to\infty}{\overset{\rm stably}{\longrightarrow}}\frac{\sigma_{1/4}}{\sqrt{2}}\int_0^1 \partial_{12} 
f(B^{(1)}_s,B^{(2)}_s)dW_s\notag\\
&&\hskip7cm+\frac14\int_0^1 \partial_{1122}f(B^{(1)}_s,B^{(2)}_s)ds\label{sn5}\\
S_n^{(6)}&:=&\sum_{k=0}^{n-1} \partial_{112} f(B^{(1)}_{k/n},B^{(2)}_{k/n}) 
\big(\Delta B^{(1)}_{k/n}\big)^2\,\Delta B^{(2)}_{k/n}
\underset{n\to\infty}{\overset{L^2}{\longrightarrow}}-\frac12\int_0^1 \partial_{1122} f(B^{(1)}_s,B^{(2)}_s)ds\label{sn6}\\
S_n^{(7)}&:=&\sum_{k=0}^{n-1} \partial_{122} f(B^{(1)}_{k/n},B^{(2)}_{k/n}) 
\Delta B^{(1)}_{k/n}\,\big(\Delta B^{(2)}_{k/n}\big)^2
\underset{n\to\infty}{\overset{L^2}{\longrightarrow}}-\frac12\int_0^1 \partial_{1122} f(B^{(1)}_s,B^{(2)}_s)ds\label{sn7}
\end{eqnarray}
\begin{eqnarray}
S_n^{(8)}&:=&\sum_{k=0}^{n-1} \partial_{1122} f(B^{(1)}_{k/n},B^{(2)}_{k/n}) 
\big(\Delta B^{(1)}_{k/n}\big)^2\,\big(\Delta B^{(2)}_{k/n}\big)^2
\underset{n\to\infty}{\overset{L^2}{\longrightarrow}}\int_0^1 \partial_{1122} f(B^{(1)}_s,B^{(2)}_s)ds\label{sn8}\\
S_n^{(9)}&:=&\sum_{k=0}^{n-1} \partial_{1112} f(B^{(1)}_{k/n},B^{(2)}_{k/n}) 
\big(\Delta B^{(1)}_{k/n}\big)^3\,\Delta B^{(2)}_{k/n}
\underset{n\to\infty}{\overset{\rm Prob}{\longrightarrow}}0\label{sn9}\\
S_n^{(10)}&:=&\sum_{k=0}^{n-1} \partial_{1222} f(B^{(1)}_{k/n},B^{(2)}_{k/n}) 
\Delta B^{(1)}_{k/n}\,\big(\Delta B^{(2)}_{k/n}\big)^3
\underset{n\to\infty}{\overset{\rm Prob}{\longrightarrow}}0\label{sn10}.
\end{eqnarray} 
Note that the term corresponding to the remainder in (\ref{ito}) converges in probability to zero 
due to the fact that $B$ has a finite quartic variation.\\

{\it Proof of (\ref{sn1}), (\ref{sn3}), (\ref{sn6}) and (\ref{sn7})}.
By Lemma \ref{lm1} with $q=3$ and $\alpha=0$, and by using the basic fact that 
\begin{equation}\label{h3}
\big(\Delta B^{(2)}_{k/n}\big)^3 = n^{-3/4}\,H_3(n^{1/4}\Delta B^{(2)}_{k/n}) + \frac{3}{\sqrt{n}}\,\Delta B^{(2)}_{k/n},
\end{equation}
we immediately see that (\ref{sn3}) is
a consequence of the following convergence:
\begin{equation}\label{sn3bis}
\frac{1}{\sqrt{n}}\sum_{k=0}^{n-1} \partial_{222} f(B^{(1)}_{k/n},B^{(2)}_{k/n}) \Delta B^{(2)}_{k/n}
\underset{n\to\infty}{\overset{L^2}{\longrightarrow}}-\frac12\int_0^1 \partial_{2222} f(B^{(1)}_s,B^{(2)}_s)ds.
\end{equation}
So, let us prove (\ref{sn3bis}). We have, on one hand:
\begin{eqnarray*}
&&E\left|\frac{1}{\sqrt{n}}\sum_{k=0}^{n-1} \partial_{222} f(B^{(1)}_{k/n},B^{(2)}_{k/n}) \Delta B^{(2)}_{k/n}\right|^2\\
&=&\frac1n\sum_{k,l=0}^{n-1} E\left(\partial_{222} f(B^{(1)}_{k/n},B^{(2)}_{k/n}) \,
\partial_{222} f(B^{(1)}_{l/n},B^{(2)}_{l/n}) \,\Delta B^{(2)}_{k/n}\,\Delta B^{(2)}_{l/n}\right)\\
&=&\frac1n\sum_{k,l=0}^{n-1} E\left(\partial_{222} f(B^{(1)}_{k/n},B^{(2)}_{k/n}) \,
\partial_{222} f(B^{(1)}_{l/n},B^{(2)}_{l/n})\,I^{(2)}_2(\delta_{k/n}\otimes\delta_{l/n})\right)\\
&+&\frac1n\sum_{k,l=0}^{n-1} E\left(\partial_{222} f(B^{(1)}_{k/n},B^{(2)}_{k/n}) \,
\partial_{222} f(B^{(1)}_{l/n},B^{(2)}_{l/n})\right)\langle\delta_{k/n},\delta_{l/n}\rangle_\HH\\
&=&\frac1n\sum_{k,l=0}^{n-1} E\left(\partial_{22222} f(B^{(1)}_{k/n},B^{(2)}_{k/n}) \,
\partial_{222} f(B^{(1)}_{l/n},B^{(2)}_{l/n})\right)
\langle\e_{k/n},\delta_{k/n}\rangle_\HH\langle\e_{k/n},\delta_{l/n}\rangle_\HH
\\
&+&\frac1n\sum_{k,l=0}^{n-1} E\left(\partial_{2222} f(B^{(1)}_{k/n},B^{(2)}_{k/n}) \,
\partial_{2222} f(B^{(1)}_{l/n},B^{(2)}_{l/n})\right)
\langle\e_{k/n},\delta_{k/n}\rangle_\HH\langle\e_{l/n},\delta_{l/n}\rangle_\HH
\\
&+&\frac1n\sum_{k,l=0}^{n-1} E\left(\partial_{2222} f(B^{(1)}_{k/n},B^{(2)}_{k/n}) \,
\partial_{2222} f(B^{(1)}_{l/n},B^{(2)}_{l/n})\right)
\langle\e_{l/n},\delta_{k/n}\rangle_\HH\langle\e_{k/n},\delta_{l/n}\rangle_\HH
\\
&+&\frac1n\sum_{k,l=0}^{n-1} E\left(\partial_{222} f(B^{(1)}_{k/n},B^{(2)}_{k/n}) \,
\partial_{22222} f(B^{(1)}_{l/n},B^{(2)}_{l/n})\right)
\langle\e_{l/n},\delta_{k/n}\rangle_\HH\langle\e_{l/n},\delta_{l/n}\rangle_\HH
\end{eqnarray*}
\begin{eqnarray*}
&+&\frac1n\sum_{k,l=0}^{n-1} E\left(\partial_{222} f(B^{(1)}_{k/n},B^{(2)}_{k/n}) \,
\partial_{222} f(B^{(1)}_{l/n},B^{(2)}_{l/n})\right)\langle\delta_{k/n},\delta_{l/n}\rangle_\HH\\
&=&a^{(n)}+b^{(n)}+c^{(n)}+d^{(n)}+e^{(n)}.
\end{eqnarray*}
Using Lemma \ref{lemma:Tec1} (i) and (ii), we have that 
$a^{(n)}$, $c^{(n)}$ and $d^{(n)}$ tends to zero as $n\to\infty$.
Using Lemma \ref{lemma:Tec1} (iii), we have
that $e^{(n)}$ tends to zero as $n\to\infty$. Finally, observe that
\begin{eqnarray*}
b^{(n)}&=& \frac1{4n^2}\sum_{k,l=0}^{n-1} E\left(\partial_{2222} f(B^{(1)}_{k/n},B^{(2)}_{k/n}) \,
\partial_{2222} f(B^{(1)}_{l/n},B^{(2)}_{l/n})\right)\\
&-&\frac1{2n\sqrt{n}}\sum_{k,l=0}^{n-1} E\left(\partial_{2222} f(B^{(1)}_{k/n},B^{(2)}_{k/n}) \,
\partial_{2222} f(B^{(1)}_{l/n},B^{(2)}_{l/n})\right)
\left(\langle\e_{l/n},\delta_{l/n}\rangle_\HH
+\frac{1}{2\sqrt{n}}\right)\\
&+&\frac1{n}\sum_{k,l=0}^{n-1} E\left(\partial_{2222} f(B^{(1)}_{k/n},B^{(2)}_{k/n}) \,
\partial_{2222} f(B^{(1)}_{l/n},B^{(2)}_{l/n})\right)
\left(\langle\e_{k/n},\delta_{k/n}\rangle_\HH
+\frac{1}{2\sqrt{n}}\right)\langle\e_{l/n},\delta_{l/n}\rangle_\HH.
\end{eqnarray*}
Therefore, using Lemma \ref{lemma:Tec1} (i) and (iv), we have
\begin{equation}\label{ehbaby}
E\left|\frac{1}{\sqrt{n}}\sum_{k=0}^{n-1} \partial_{222} f(B^{(1)}_{k/n},B^{(2)}_{k/n}) \Delta B^{(2)}_{k/n}\right|^2
=E\left|\frac{1}{2n}\sum_{k=0}^{n-1}\partial_{2222}  f(B^{(1)}_{k/n},B^{(2)}_{k/n}) \right|^2 + o(1).
\end{equation}

On the other hand, we have
\begin{eqnarray*}
&&E\left(\frac{1}{\sqrt{n}}\sum_{k=0}^{n-1} \partial_{222} f(B^{(1)}_{k/n},B^{(2)}_{k/n}) \Delta B^{(2)}_{k/n}
\times \frac{-1}{2n}\sum_{l=0}^{n-1} \partial_{2222} f(B^{(1)}_{l/n},B^{(2)}_{l/n}) \right)\\
&=&-\frac1{2n\sqrt{n}}\sum_{k,l=0}^{n-1} E\big(\partial_{222}f(B^{(1)}_{k/n},B^{(2)}_{k/n})
\partial_{2222}f(B^{(1)}_{l/n},B^{(2)}_{l/n})\, \Delta B^{(2)}_{k/n}\big)\\
&=&-\frac1{2n\sqrt{n}}\sum_{k,l=0}^{n-1} E\big(\partial_{2222}f(B^{(1)}_{k/n},B^{(2)}_{k/n})
\partial_{2222}f(B^{(1)}_{l/n},B^{(2)}_{l/n})\big)\langle \e_{k/n},\delta_{k/n}\rangle_\HH\\
&-&\frac1{2n\sqrt{n}}\sum_{k,l=0}^{n-1} E\big(\partial_{222}f(B^{(1)}_{k/n},B^{(2)}_{k/n})
\partial_{22222}f(B^{(1)}_{l/n},B^{(2)}_{l/n})\big)\langle \e_{l/n},\delta_{k/n}\rangle_\HH\\
&=&\frac1{4n^2}\sum_{k,l=0}^{n-1} E\big(\partial_{2222}f(B^{(1)}_{k/n},B^{(2)}_{k/n})
\partial_{2222}f(B^{(1)}_{l/n},B^{(2)}_{l/n})\big)\\
&-&\frac1{2n\sqrt{n}}\sum_{k,l=0}^{n-1} E\big(\partial_{2222}f(B^{(1)}_{k/n},B^{(2)}_{k/n})
\partial_{2222}f(B^{(1)}_{l/n},B^{(2)}_{l/n})\big)\left(\langle \e_{k/n},\delta_{k/n}\rangle_\HH+\frac1{2\sqrt{n}}\right)\\
&-&\frac1{2n\sqrt{n}}\sum_{k,l=0}^{n-1} E\big(\partial_{222}f(B^{(1)}_{k/n},B^{(2)}_{k/n})
\partial_{22222}f(B^{(1)}_{l/n},B^{(2)}_{l/n})\big)\langle \e_{l/n},\delta_{k/n}\rangle_\HH
\end{eqnarray*}
We immediately have that the second (see Lemma \ref{lemma:Tec1} (iv)) and the third
(see Lemma \ref{lemma:Tec1} (ii)) terms in the previous expression tends to
zero as $n\to\infty$. That is
\begin{eqnarray}
&&E\left(\frac{1}{\sqrt{n}}\sum_{k=0}^{n-1} \partial_{222} f(B^{(1)}_{k/n},B^{(2)}_{k/n}) \Delta B^{(2)}_{k/n}
\times \frac{-1}{2n}\sum_{l=0}^{n-1} \partial_{2222} f(B^{(1)}_{l/n},B^{(2)}_{l/n}) \right)\notag\\
&=&E\left|\frac{1}{2n}\sum_{k=0}^{n-1}\partial_{2222}  f(B^{(1)}_{k/n},B^{(2)}_{k/n}) \right|^2 + o(1).
\label{ehehbaby}
\end{eqnarray}

We have proved, see (\ref{ehbaby}) and (\ref{ehehbaby}), that
$$
E\left|\frac{1}{\sqrt{n}}\sum_{k=0}^{n-1} \partial_{222} f(B^{(1)}_{k/n},B^{(2)}_{k/n}) \Delta B^{(2)}_{k/n}
+\frac{1}{2n}\sum_{k=0}^{n-1} \partial_{2222} f(B^{(1)}_{k/n},B^{(2)}_{k/n}) \right|^2
\underset{n\to\infty}{\longrightarrow} 0.
$$
This implies (\ref{sn3bis}).

The proof of (\ref{sn1}) follows directly from (\ref{sn3})
by exchanging the roles played by $B^{(1)}$ and $B^{(2)}$.
On the other hand, by combining Lemma \ref{lm1} with the following basic identity: 
$$\big(\Delta B^{(2)}_{k/n}\big)^2 = \frac{1}{\sqrt{n}}\,H_2(n^{1/4}\Delta B^{(2)}_{k/n}) + \frac{1}{\sqrt{n}},$$ 
we see that (\ref{sn7}) is also a direct consequence of (\ref{sn3bis}).
Finally, (\ref{sn6}) is obtained from (\ref{sn7}) by exchanging the roles played by $B^{(1)}$ and $B^{(2)}$.\\

{\it Proof of (\ref{sn2}), (\ref{sn4}) and (\ref{sn8})}.
By combining Lemma \ref{lm1} with the identity
$$\big(\Delta B^{(1)}_{k/n}\big)^4 = \frac1n\,H_4(n^{1/4}\Delta B^{(1)}_{k/n}) +
\frac{6}{n}\,H_2(n^{1/4}\Delta B^{(1)}_{k/n}) +\frac3n,$$
we see that (\ref{sn4}) is easily obtained through a Riemann sum argument.
We can use the same arguments in order to prove (\ref{sn2}). Finally, to obtain (\ref{sn8}), it suffices to combine Lemma \ref{lm1}
with the identity
$$
\big(\Delta B^{(1)}_{k/n}\big)^2 \big(\Delta B^{(2)}_{k/n}\big)^2
= \frac{1}{n} + \frac{1}{\sqrt{n}}\,\big(\Delta B^{(1)}_{k/n}\big)^2\,H_2(n^{1/4}\Delta B^{(2)}_{k/n})
+\frac{1}{n}\,H_2(n^{1/4}\Delta B^{(1)}_{k/n}).
$$

{\it Proof of (\ref{sn9}) and (\ref{sn10})}. We only prove (\ref{sn10}), the proof of (\ref{sn9})
being obtained from (\ref{sn10}) by exchanging the roles played by $B^{(1)}$ and $B^{(2)}$.
By combining (\ref{h3}) with Lemma \ref{lm1}, it
suffices to prove that
$$
\frac1{\sqrt{n}}\sum_{k=0}^{n-1}\partial_{1222}f(B^{(1)}_{k/n},B^{(2)}_{k/n})\,\Delta B^{(1)}_{k/n}\,
\Delta B^{(2)}_{k/n}\underset{n\to\infty}{\overset{\rm Prob}{\longrightarrow}} 0.
$$
But this last convergence follows directly from Lemma \ref{law}. Therefore, the proof of (\ref{sn10})
is done.\\

{\it Proof of (\ref{sn5})}. We combine Proposition \ref{law} with the idea developed in the third comment
that
we have addressed just after the statement of Theorem \ref{main-thm}.
Indeed, we have
\begin{eqnarray*}
&&\sum_{k=0}^{n-1}\partial_{12}f\big(B^{(1)}_{k/n},B^{(2)}_{k/n}\big)\Delta B^{(1)}_{k/n}\Delta B^{(2)}_{k/n}\\
&=&
\frac{1}{2\sqrt{n}}\sum_{k=0}^{n-1}\partial_{12}f
\big(
\frac{\beta_{k/n}+\widetilde{\beta}_{k/n}}{\sqrt 2},
\frac{\beta_{k/n}-\widetilde{\beta}_{k/n}}{\sqrt 2}
\big)
\left(
\sqrt{n}\big(\Delta\beta_{k/n}\big)^2- 1
\right)\\
&&-
\frac{1}{2\sqrt{n}}\sum_{k=0}^{n-1}\partial_{12}f
\big(
\frac{\beta_{k/n}+\widetilde{\beta}_{k/n}}{\sqrt 2},
\frac{\beta_{k/n}-\widetilde{\beta}_{k/n}}{\sqrt 2}
\big)
\left(
\sqrt{n}\big(\Delta\widetilde{\beta}_{k/n}\big)^2- 1
\right)
\end{eqnarray*}
for $\beta=(B^{(1)}+B^{(2)})/\sqrt{2}$ and $\widetilde{\beta}=(B^{(1)}-B^{(2)})/\sqrt{2}$.
Note that $(\beta,\widetilde{\beta})$ is also a 2D fractional Brownian motion of Hurst index $1/4$.
Hence, using Proposition \ref{law} with 
$
g(x,y)=
\widetilde{g}(x,y)=
f\left(
\frac{x+y}{\sqrt{2}}
,
\frac{x-y}{\sqrt{2}}
\right)
$, we get
\begin{eqnarray*}
&&\sum_{k=0}^{n-1}\partial_{12}f\big(B^{(1)}_{k/n},B^{(2)}_{k/n}\big)\Delta B^{(1)}_{k/n}\Delta B^{(2)}_{k/n}\\
&\underset{n\to\infty}{\overset{\rm stably}{\longrightarrow}}&
\frac{\sigma_{1/4}}{2}\int_0^1 \partial_{12}f(B^{(1)}_s,B^{(2)}_s)d(W-\widetilde{W})_s
+\frac14\int_0^1 \partial_{1122} f(B^{(1)}_s,B^{(2)}_s)ds\\
&\overset{\rm Law}{=}&\frac{\sigma_{1/4}}{\sqrt{2}}\int_0^1 \partial_{12}f(B^{(1)}_s,B^{(2)}_s)dW_s
+\frac14\int_0^1 \partial_{1122} f(B^{(1)}_s,B^{(2)}_s)ds,
\end{eqnarray*}
for $(W,\widetilde{W})$ a 2D standard Brownian motion independent of $(\beta,\widetilde{\beta})$.
The proof of (\ref{sn5}) is done.\\
\\
{\bf Proof of the first point} (case $H>1/4$). The proof can be done by following exactly the same strategy than in the step above.
The only difference is that, using a version of Lemma \ref{lm1} together with computations similar to that 
allowing to obtain (\ref{sn3bis}), the limits in (\ref{sn1})--(\ref{sn8}) are, here, {\sl all equal to zero}
(for the sake of simplicity, the technical details are left to the reader). Therefore, 
we can deduce (\ref{h>1/4}) by using (\ref{ito}).
\fin

\end{document}